\documentclass[11pt,twoside,onecolumn]{article}

\usepackage{authblk} 

\usepackage{graphicx}
\graphicspath{ {Vis/} }
\usepackage[hang, small,labelfont=bf,up,textfont=it,up]{caption}
\usepackage{subcaption,dsfont}

\usepackage[usenames, dvipsnames]{color}

\usepackage{indentfirst} 

\usepackage[sc]{mathpazo} 
\usepackage[T1]{fontenc} 
\usepackage{microtype} 

\usepackage[english]{babel} 

\usepackage[hmarginratio=1:1,top=32mm,columnsep=20pt]{geometry} 
\usepackage{booktabs} 

\usepackage{lettrine} 

\usepackage{enumitem} 
\setlist[itemize]{noitemsep} 

\usepackage{titlesec} 

\usepackage{fancyhdr} 
\usepackage{hyperref} 

\usepackage{mathtools} 

\usepackage{amsmath,calc,amsthm}

\DeclareMathOperator*{\argmin}{argmin}

\usepackage{amsfonts}
\usepackage{amsthm}

\theoremstyle{definition}

\usepackage{multicol}
\usepackage{multirow}
\usepackage{amssymb}
\usepackage{array}
\usepackage{tabu}
\usepackage{algpseudocode}
\usepackage[numbers]{natbib}
\usepackage{url}

\theoremstyle{definition}

\theoremstyle{plain}

\usepackage{xcolor}

\usepackage{makecell}
\usepackage{scrextend} 

\usepackage{algorithm}
\usepackage{nomencl}
\usepackage{etoolbox}
\makenomenclature

\makeatletter
\patchcmd{\thenomenclature}{\section*}{\section}{}{}
\makeatother

\begin{document}

\title{Integrating Condition-Based Maintenance\\into Dynamic Spare Parts Management} 

\author[1]{D. Usanov\thanks{Corresponding author. Tel.:+31(0)20 592 4168.\\
\textit{E-mail addresses}: usanov@cwi.nl (D. Usanov), p.m.van.de.ven@cwi.nl (P.M. van de Ven), r.d.van.der.mei@cwi.nl (R.D. van der Mei).}}
\author[1]{P.M. van de Ven}
\author[1, 2]{R.D. van der Mei}
\affil[1]{Centrum Wiskunde \& Informatica, Science Park 123,
1098 XG Amsterdam, The Netherlands}
\affil[2]{Vrije Universiteit Amsterdam, De Boelelaan 1081a,
1081 HV Amsterdam, The Netherlands }

\date{} 


\maketitle
\newpage
\begin{abstract}

In this paper we introduce a new model where the concept of condition-based maintenance
is combined in a network setting with dynamic spare parts management. The
model facilitates both preventive and corrective maintenance of geographically distributed
capital goods as well as relocation of spare parts between different warehouses based on
the availability of stock and the condition of all capital good installations. We formulate
the problem as a Markov decision process, with the degradation process explicitly incorporated
into the model. Numerical experiments show that that significant cost savings can
be achieved when condition monitoring is used for preventive maintenance in a service
network for capital goods.

\end{abstract}



\section{Introduction}

Capital goods, such as MRI scanners, lithography machines, aircraft, or wind turbines are subject to deterioration and require maintenance over their lifetime. Continuous operation of such assets is crucial, as failures can have significant negative effects. Advancements in condition monitoring techniques facilitate tracking the degradation process of capital goods in real time. This creates a tremendous potential for implementing preventive maintenance policies that use sensor data to indicate which spare parts should be replaced before a breakdown happens. This type of preventive maintenance is called \textit{condition-based maintenance} (CBM), and it can be extremely useful to mitigate the risks related to downtime of capital goods.

The existing research on CBM is focused on optimizing control limits and/or maintenance intervals for one single- or multi-component machine. However, these works typically ignore the fact that these machines are parts of a network comprising many machines distributed across different geographical locations, and the spare parts supply distributed over a number of stock points to ensure short response times to failures. CBM policies developed for a single machine do not take into account the location and availability of stock, as well as the location and condition of other machines that might require maintenance. The research literature on dynamic spare parts management in a service network is primarily focused on optimal corrective maintenance and relocation of spare parts, where the failures are typically assumed to follow a Poisson process and cannot be predicted.

\begin{figure}
\begin{subfigure}{.4\textwidth}
\centering
\includegraphics[width=\textwidth]{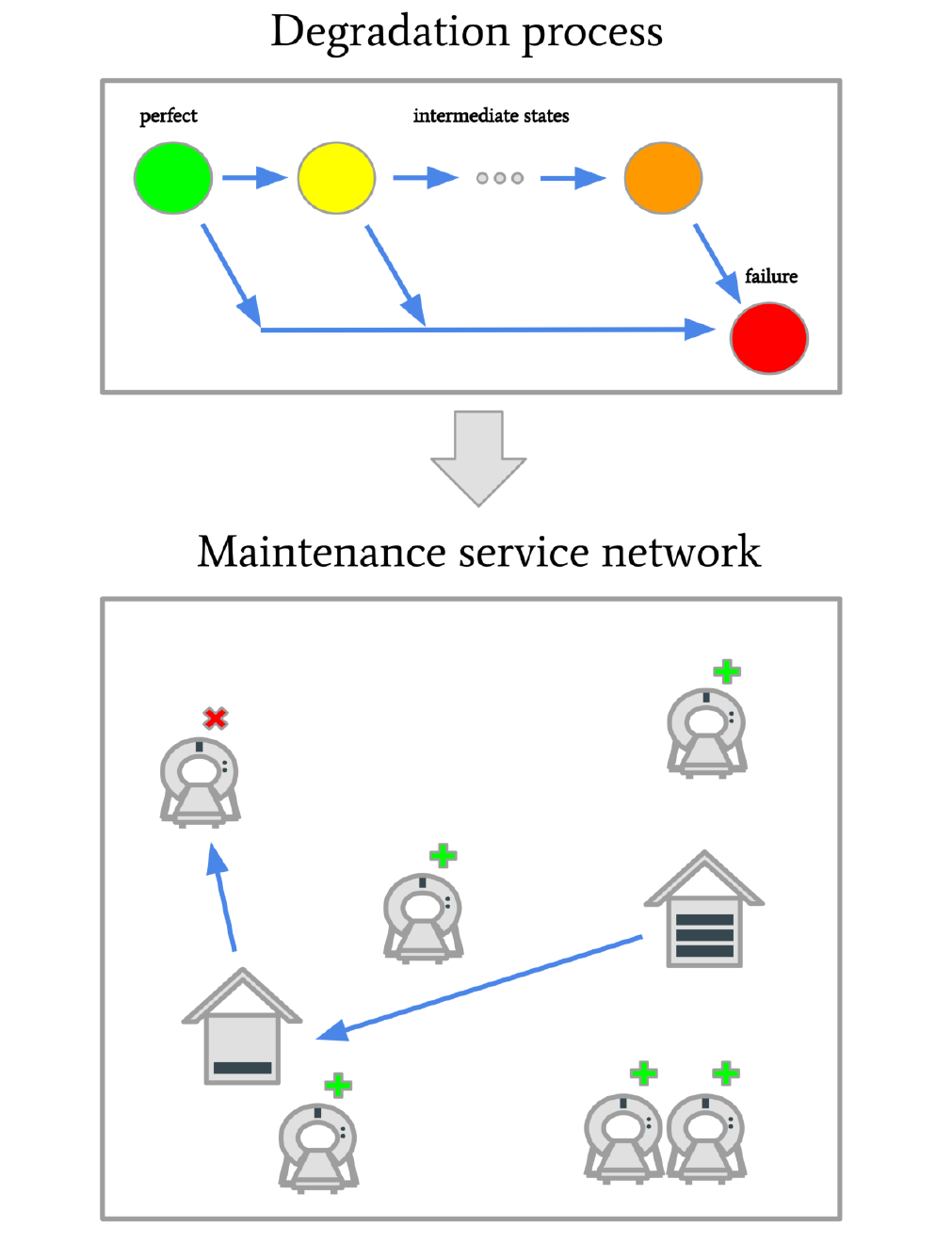}
\caption{}
\label{fig:cbm_intro_1}
\end{subfigure}\hfill
\begin{subfigure}{.55\textwidth}
\centering
\includegraphics[width=\textwidth]{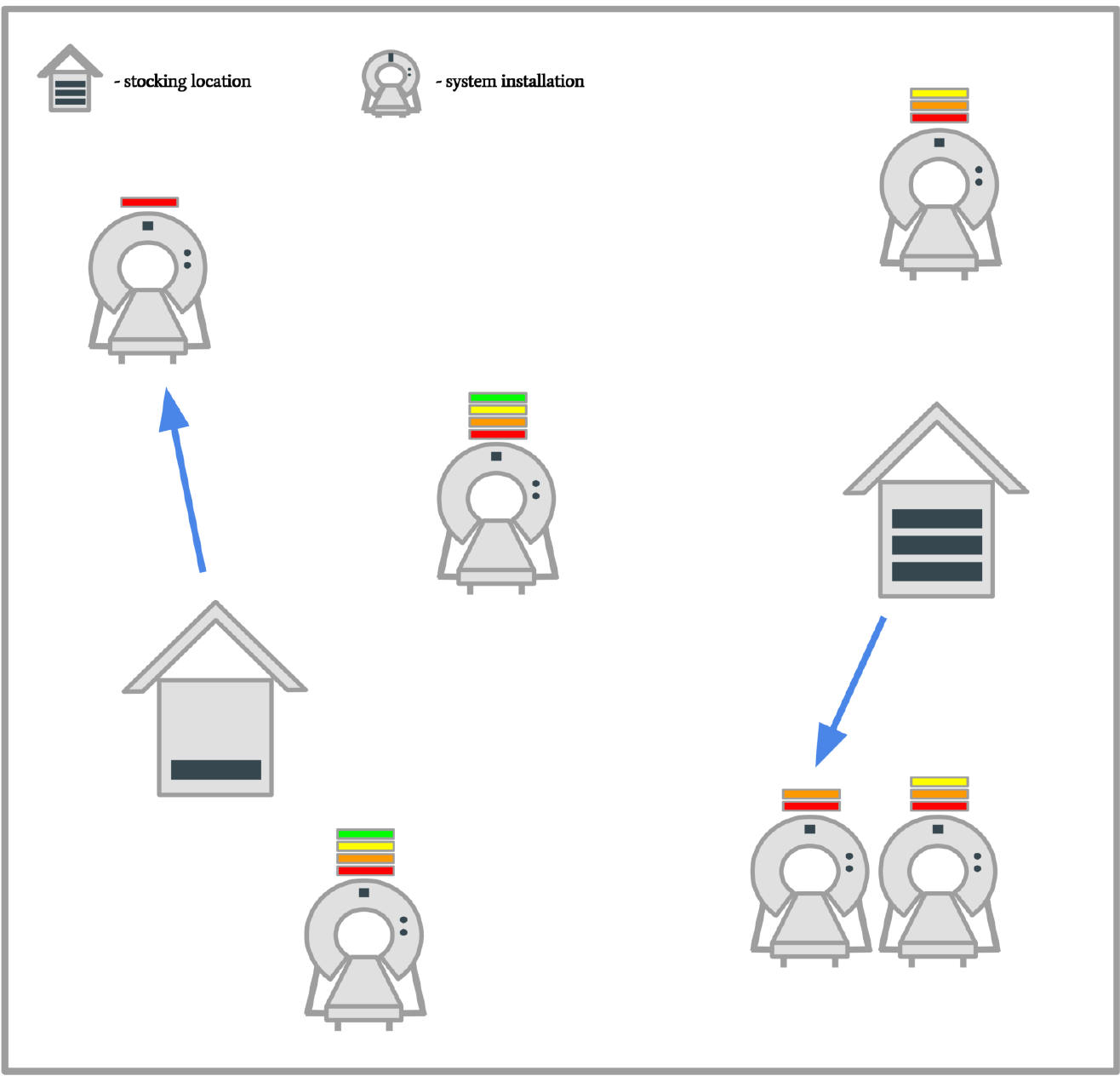}
\caption{}
\label{fig:cbm_intro_2}
\end{subfigure}
	\caption{Integrating CBM into dynamic spare parts management}
    \label{fig:cbm_intro}
\end{figure}

In this paper we integrate the CBM concept into a network setting with dynamic spare parts management. This allows us to relocate spare parts between stock points and perform proactive maintenance of machines based on stock levels and the condition of all machines. Figure~\ref{fig:cbm_intro} illustrates how incorporating CBM changes the complexity and the dynamics of the service network. The state space increases, as each of the machines has more than just two possible condition states (perfect and failure), as it is typically assumed in the research literature on dynamic spare parts management. However, this provides more information about the overall state of the network, and therefore, enables more educated decision making. For instance, instead of relocating a spare part upon failure of one of the machines, it might be better to preventively repair another machine that is close to failure. Using the information obtained through condition monitoring can improve both maintenance and relocation activities, and boost the maximum performance of the service network.

We consider a single-component machine and assume a Markovian degradation process, where a machine moves through a sequence of intermediate states before it reaches the failure state. This is a common assumption in research literature (see for example~\cite{kharoufeh2010,jiang2013}). A number of such machines are spread across a service region, and we optimize corrective and preventive maintenance actions, as well as proactive relocation of spare parts between stock points. We show that, by introducing condition monitoring into a network setting, significant improvements can be achieved in reducing total expected costs, independent of the cost structure.

To summarize, in this paper we make the following contributions:
\begin{enumerate}
\item We propose a new model, where condition-based maintenance is integrated into the dynamic spare-parts management;
\item We formulate the problem as a Markov Decision Process (MDP), and conduct numerical experiments showing that incorporating CBM can significantly reduce the maintenance costs.
\end{enumerate}

The remainder of this paper is organized as follows. In Section~\ref{sec:lit}, we provide brief review of the relevant literature. Section~\ref{sec:model}, presents the model and the MDP formulation. In Section~\ref{sec:num} we numerically evaluate the impact of CBM the optimal policy performance. Conclusions and suggestions for further research are made in Section~\ref{sec:conclusion}.

\section{Literature Review}\label{sec:lit}

The research most related to our work comes from the two streams of literature: dynamic dispatching and relocation of spare parts, and condition-based maintenance. Although attention to both topics seems to increase in recent years, to the best of our knowledge, there has been no research on combining these two concepts. Below we make a brief overview of the latest work in the two domains.

Dispatching and relocation of spare parts considers operational level decision making in maintenance service networks. Proactive and reactive allocation of stock in spare parts networks is commonly referred to in research literature as \textit{lateral transshipment}. For a comprehensive overview of the research done on lateral transshipments we refer to~\cite{wong2006, paterson2011}. In~\cite{wijk2009} the authors consider an inventory model with fixed inventory level and two warehouses, both facing Poisson demand for spare parts. They provide an exact analysis of the model that derives an optimal policy for allocation of demand to warehouses. In~\cite{tiemessen2013} a dynamic demand allocation rule is developed that is scalable for spare parts networks of realistic size. The authors show that significant cost savings can be achieved with their approach compared to the static allocation rule commonly used in practice, where the closest warehouse is used to fulfill the demand. An interesting contribution is made in~\cite{paterson2012}, where the authors consider relocating additional stock when satisfying real-time demand. Recent developments include~\cite{feng2017, meissner2018}, where proactive relocation of stock is studied along with the reactive policies.

The recent works on condition-based maintenance include~\cite{zhu2015, peng2016, zhu2017, drent2019}. In~\cite{zhu2015} the authors study a multi-component system where each component follows a stochastic degradation process according to the so-called random coefficient model. A joint maintenance of components whose condition falls below a control limit is performed at scheduled downs. Those control limits together with the maintenance interval are subject to optimization. In~\cite{zhu2017} the authors consider condition-based maintenance of a single component that is a part of a complex system. This component follows a stochastic degradation process for which the authors use the random coefficient model and Gamma process. A control limit policy is analyzed, with maintenance actions taken at scheduled or unscheduled downs related to other components of the system. A single component model with stochastic degradation process is also studied in~\cite{peng2016}. The authors consider an application for a manufacturing system and look into joint optimization of control limit and production quantity of the lot-sizing policy. Another way to model the degradation process is using Markov process with discrete states. In~\cite{drent2019} the authors consider a single-component system that follows Markovian degradation process with two intermediate states and study the optimal control limit policy of such system.

Our work is different from these streams of research literature, as we consider a generic model where \textit{both reactive and proactive allocation of stock} is allowed in real-time, along with condition-based maintenance.

\section{Model}\label{sec:model}

We consider a network of identical single-component machines supported by a set of local warehouses and a central warehouse. The state of each machine is completely and continuously observable. Replacement of components happens either \textit{preventively} or \textit{correctively} upon a failure. Let $\mathcal{I} = \{0, 1, ..., I\}$ be a set of warehouses, with $i = 0$ - central warehouse with ample capacity. Let $\mathcal{J} = \{1, ..., J\}$ be a set of machines.

We assume that the lifetime of a machine is Cox distributed with $N>0$ phases~\cite{koole2004}. We choose Cox distribution because it allows to approximate any random variable with positive support. Denote the condition of a machine $j \in \mathcal{J}$ by $C_j \in \{0, 1, ..., N\}$, where 0 corresponds to failure and $N$ to the perfect condition. A machine stays in each state $n \in \{1, ..., N\}$ for an exponential amount of time with parameter $\mu_n$, then it moves either to the 'failure' state 0 with probability $\alpha_n < 1$ or to the state $n-1$ with probability $1-\alpha_n$. Note that from state 1 the machine moves to state 0 with probability 1, so $\alpha_1 = 1$. Upon breakdown, a spare part is dispatched to that machine from either a local or the central warehouse. The downtime of a machine due to corrective maintenance includes traveling time and repair time, and we assume it is exponentially distributed  with parameter $\mu_0$. A machine can be preventively repaired at any point in time. We assume that there is no downtime of a machine in this case, so a spare part is replaced instantly. This is a common assumption in literature, as preventive maintenance is typically easier than corrective, and the corresponding downtime does not include delivery of a spare part. After a spare part is replaced, either correctively or preventively, a machine moves back into the `perfect' state $N$.

If a failure occurs, a spare part is to be dispatched immediately either from a local warehouse, or from the central warehouse. 
Once a spare part is dispatched from a local warehouse, a replenishment order is placed. The replenishment lead time is exponentially distributed with parameter $\gamma$. If a spare part was dispatched from a local warehouse $i$, a relocation of a single spare part from one of the other warehouses to the warehouse $i$ is allowed. We assume that such relocations happen instantaneously. This is a reasonable assumption, as the traveling times between warehouses are typically low compared to the average time between consecutive failures of capital goods.

If the condition of a machine degrades but it is still not in a failure state, we consider two types of actions. We may decide to repair this machine preventively, and do up to one relocation the local warehouse from which a spare part is dispatched. Alternatively, we may decide to not repair the machine, and do up to one relocation between any two local warehouses. These relocations are intended to better distribute stock across the network, and reduce future downtime.

Let $\boldsymbol C = (C_1,...,C_J)$, $\boldsymbol F = (F_1,...,F_I)$, $\boldsymbol P = (P_1,...,P_I)$, where $F_i \geq 0$ and $P_i \geq 0$ denote the stock level and the pipeline stock (replenishment orders) at warehouse $i$, respectively. Denote by $K = \sum_{i \in \mathcal{I}}\left(F_i(t)+P_i(t)\right)$ the aggregate inventory level. Note that $K$ always remains constant, as each time a spare part is dispatched, a new one is ordered immediately and added to the corresponding pipeline stock. Let $\boldsymbol X = (\boldsymbol F, \boldsymbol P, \boldsymbol C, j)$ with $j\in \{0, 1\dots, J\}$, denote the state of the system immediately after the condition of machine $j \in \mathcal{J}$ changes or a replenishment order arrives ($j=0$). Let also $\boldsymbol a(\boldsymbol X) = (x, y, z)$ represent the action in state $\boldsymbol X$. Here, $x \in \{-1, 0, 1,...,I\}$ indicates the warehouse from which a spare part is to be dispatched, $y \in \{-1, 1,...,I\}$ indicates the warehouse from which a spare part is to be relocated, and $z \in \{-1, 1,...,I\}$ the warehouse to which a relocated spare part should be placed. The value $x = -1$ corresponds to the case when no dispatching is made, $y = z = -1$ to the decision not to relocate a spare part.

\subsection{Actions}\label{sec:model:actions}

We consider two types of actions. The first type includes both dispatching a spare part and relocating another spare part to the warehouse from which the dispatching was made. Relocating a spare part is not necessary and is only considered if a part was dispatched from a local warehouse. This type of action can be made either correctively upon a failure, in which case repair is required, or preventively when a machine's condition degrades but the machine is still functioning. Let $\mathcal{W}(\pmb{X}) \subseteq \mathcal{I}$ denote the set of local warehouses that have at least one spare part in stock in state $\pmb{X}$. For a state $\pmb{X} = (\boldsymbol F, \boldsymbol P, \boldsymbol C, j)$, $C_j < N$, the type-1 action space is defined as
\begin{flalign}\label{eq:mdp_actions1}
\mathcal{A}_1(\boldsymbol X) = &\{(x, y, z)|x \in \mathcal{W}(\pmb{X}), \ y \in \mathcal{W}(\pmb{X})\setminus \{x\}, \ z = x\} \\
\nonumber &\cup \{(x, y, z)|x \in \mathcal{W}(\pmb{X})\cup \{0\}, \ y = -1, \ z = -1\}.
\end{flalign}%
The second type of action includes only relocations. These are allowed upon a change of a machine state that does not result in a failure. In this case, a single relocation is allowed between any pair of local warehouses as long as the origin warehouse has a spare part available. For a state $\pmb{X} = (\boldsymbol F, \boldsymbol P, \boldsymbol C, j)$, $C_j > 0$, the type-2 action space is defined as
\begin{flalign}\label{eq:mdp_actions2}
\mathcal{A}_2(\boldsymbol X) = &\{(x, y, z)|x = -1, \ y \in \mathcal{W}(\pmb{X}), \ z \in \mathcal{I}\setminus \{y\}\} \cup \{(-1, -1, -1)\}.
\end{flalign}%
Thus, the total action space is $\mathcal{A}(\boldsymbol X) = \mathcal{A}_1(\boldsymbol X) \cup \mathcal{A}_2(\boldsymbol X)$.

We denote by $\boldsymbol X^{\boldsymbol a}(t)$ the state of the system at time $t$ under decision rule $\boldsymbol a$. The process $\{\boldsymbol X^{\boldsymbol a}(t)\}_{t \geq 0}$ is a continuous-time Markov process, with the following state space:
\begin{flalign*}
\mathcal{S} = \big\{(\boldsymbol F, \boldsymbol P, \boldsymbol C, j)| & F_i \geq 0, P_i \geq 0 \ \forall i \in \mathcal{I}; \sum_{i \in \mathcal{I}}\left(F_i+P_i\right) = K; \\
& C_j \in \{0,1,...,N\} \ \forall j \in \mathcal{J}; j \in \{0, 1,...,J\}\big\}.
\end{flalign*}%

\subsection{Transitions}\label{sec:model:transitions}
\noindent
In this subsection we define transition rates between states corresponding to decision epochs. We consider the following four types of events when the state of the system changes:
\begin{enumerate}
\item arrival of a replenishment order at local warehouse $i$;
\item failure of machine $j$;
\item degradation of machine $j$ that is not a failure;
\item repair of machine $j$.
\end{enumerate}

Given the machines' condition vector $\boldsymbol C$, let $\boldsymbol C^{j,n}$ denote a vector that is obtained from $\boldsymbol C$ by setting its $j$-th component to $n$, so $C^{j,n}_k = C_k \ \forall k \neq j$ and $C^{j,n}_j = n$. We denote as $\boldsymbol e_k$ the vector of length $I$ with the $k$-th element equal to 1 and all other elements equal to 0 for $k \in \{1,...I\}$. Finally, $\boldsymbol e_0$ and $\boldsymbol e_{-1}$ both denote a zero vector of length $I$.

Assume that the system is in state $\boldsymbol X = (\boldsymbol F, \boldsymbol P, \boldsymbol C, j)$ immediately after an event occurred at machine $j$ and before an action is taken. Remember that we set $j=0$ if the last even is an arrival of a replenishment order. We want to define the transition rate $q^{\pmb{a}}\big((\pmb{F}, \pmb{P}, \pmb{C}, j), (\pmb{F}'+, \pmb{P}', \pmb{C}', j')\big)$ from each state $(\pmb{F}, \pmb{P}, \pmb{C}, j)$ to each possible next state $(\pmb{F}'+, \pmb{P}', \pmb{C}', j')$ that is determined by an action $\pmb{a}$ and the next event $j'$. The following types of transitions are possible in our model.

\noindent
\textbf{Type 1.} \underline{\textit{Last event:}} replenishment at warehouse $i \in \{1,..,I\}$; \underline{\textit{action:}} $\pmb{a} = (-1, -1, -1)$.\\
The last event is a replenishment at a local warehouse $i$. In that case only one action is possible, and that is to do nothing. The state of the system immediately after the action is taken is $\pmb{X} = (\pmb{F}, \pmb{P}, \pmb{C}, 0)$. The transition rates then depend on the next state, and are defined as follows.
\begin{enumerate}
\item \underline{\textit{Next event:}} replenishment at $k$, $P_k>0$. The corresponding transition rate is
$$q^{\pmb{a}}\big((\pmb{F}, \pmb{P}, \pmb{C}, 0), (\pmb{F}+\pmb{e}_k, \pmb{P}-\pmb{e}_k, \pmb{C}, 0)\big) = P_k\gamma.$$
\item \underline{\textit{Next event:}} failure of machine $j$, $C_j>0$. The corresponding transition rate is
$$q^{\pmb{a}}\big((\pmb{F}, \pmb{P}, \pmb{C}, 0), (\pmb{F}, \pmb{P}, \pmb{C}^{j, 0}, j)\big) = \alpha_{C_j}\mu_{C_j}.$$
\item \underline{\textit{Next event:}} degradation of machine $j$, $C_j>1$. The corresponding transition rate is
$$q^{\pmb{a}}\big((\pmb{F}, \pmb{P}, \pmb{C}, 0), (\pmb{F}, \pmb{P}, \pmb{C}^{j, C_j-1}, j)\big) = (1-\alpha_{C_j})\mu_{C_j}.$$
\item \underline{\textit{Next event:}} repair of machine $j$, $C_j=0$. The corresponding transition rate is
$$q^{\pmb{a}}\big((\pmb{F}, \pmb{P}, \pmb{C}, 0), (\pmb{F}+\pmb{e}_k, \pmb{P}-\pmb{e}_k, \pmb{C}^{j, N}, j)\big) = \mu_0.$$
\end{enumerate}

\noindent
\textbf{Type 2.1} \underline{\textit{Last event:}} failure of machine $j \in \mathcal{J}$; \underline{\textit{action:}} $\pmb{a}=(x, y, z) \in \mathcal{A}_1(\pmb{F}, \pmb{P}, \pmb{C}, j)$.\\
The last event is a failure of machine $j$, and the action is to dispatch a spare part from a warehouse $x$. The action may also include a relocation of one spare part from warehouse $y$ to warehouse $z = x$. The state of the system immediately after the action is taken is $\pmb{X} = (\pmb{F}-\pmb{e}_x\mathbb{1}\{y=-1\}-\pmb{e}_y, \pmb{P}+\pmb{e}_x, \pmb{C}, j)$. The transition rates are defined as follows.
\begin{enumerate}
\item \underline{\textit{Next event:}} replenishment at $k$, $(\pmb{P}+\pmb{e}_x)_k>0$. The corresponding transition rate is
$$q^{\pmb{a}}\big((\pmb{F}, \pmb{P}, \pmb{C}, j), (\pmb{F}-\pmb{e}_x\mathbb{1}\{y=-1\}-\pmb{e}_y+\pmb{e}_k, \pmb{P}+\pmb{e}_x-\pmb{e}_k, \pmb{C}, 0)\big) = (\pmb{P}+\pmb{e}_x)_k\gamma.$$
\item \underline{\textit{Next event:}} failure of machine $l$, $C_l>0$. The corresponding transition rate is
$$q^{\pmb{a}}\big((\pmb{F}, \pmb{P}, \pmb{C}, j), (\pmb{F}-\pmb{e}_x\mathbb{1}\{y=-1\}-\pmb{e}_y, \pmb{P}+\pmb{e}_x, \pmb{C}^{l, 0}, l)\big) = \alpha_{C_l}\mu_{C_l}.$$
\item \underline{\textit{Next event:}} degradation of machine $l$, $C_l>1$. The corresponding transition rate is
$$q^{\pmb{a}}\big((\pmb{F}, \pmb{P}, \pmb{C}, j), (\pmb{F}-\pmb{e}_x\mathbb{1}\{y=-1\}-\pmb{e}_y, \pmb{P}+\pmb{e}_x, \pmb{C}^{l, C_l-1}, l)\big) = (1-\alpha_{C_l})\mu_{C_l}.$$
\item \underline{\textit{Next event:}} repair of machine $l$, $C_l=0$. The corresponding transition rate is
$$q^{\pmb{a}}\big((\pmb{F}, \pmb{P}, \pmb{C}, j), (\pmb{F}-\pmb{e}_x\mathbb{1}\{y=-1\}-\pmb{e}_y, \pmb{P}+\pmb{e}_x, \pmb{C}^{l, N}, l)\big) = \mu_0.$$
\end{enumerate}

\noindent
\textbf{Type 2.2} \underline{\textit{Last event:}} degradation of machine $j \in \mathcal{J}$; \underline{\textit{action:}} $\pmb{a}=(x, y, z) \in \mathcal{A}_1(\pmb{F}, \pmb{P}, \pmb{C}, j)$.\\
The last event is a degradation of machine $j$ that is not a failure. A spare part from a warehouse $x$ is dispatched for preventive maintenance of the machine. The action may also include a relocation of one spare part from warehouse $y$ to warehouse $z = x$. As preventive maintenance assumed to be done instantaneously, the state of the system immediately after the action is taken is $\pmb{X} = (\pmb{F}-\pmb{e}_x\mathbb{1}\{y=-1\}-\pmb{e}_y, \pmb{P}+\pmb{e}_x, \pmb{C}^{j, N})$. The transition rates are defined as follows.
\begin{enumerate}
\item \underline{\textit{Next event:}} replenishment at $k$, $(\pmb{P}+\pmb{e}_x)_k>0$. The corresponding transition rate is
$$q^{\pmb{a}}\big((\pmb{F}, \pmb{P}, \pmb{C}, j), (\pmb{F}-\pmb{e}_x\mathbb{1}\{y=-1\}-\pmb{e}_y+\pmb{e}_k, \pmb{P}+\pmb{e}_x-\pmb{e}_k, \pmb{C}^{j, N}, 0)\big) = (\pmb{P}+\pmb{e}_x)_k\gamma.$$
\item \underline{\textit{Next event:}} failure of machine $l$, $C_l>0$. The corresponding transition rate is
$$q^{\pmb{a}}\big((\pmb{F}, \pmb{P}, \pmb{C}, j), (\pmb{F}-\pmb{e}_x\mathbb{1}\{y=-1\}-\pmb{e}_y, \pmb{P}+\pmb{e}_x, {\pmb{C}^{j, N}}^{l, 0}, l)\big) = \alpha_{C_l}\mu_{C_l}.$$
\item \underline{\textit{Next event:}} degradation of machine $l$, $C_l>1$. The corresponding transition rate is
$$q^{\pmb{a}}\big((\pmb{F}, \pmb{P}, \pmb{C}, j), (\pmb{F}-\pmb{e}_x\mathbb{1}\{y=-1\}-\pmb{e}_y, \pmb{P}+\pmb{e}_x, {\pmb{C}^{j, N}}^{l, C_l-1}, l)\big) = (1-\alpha_{C_l})\mu_{C_l}.$$
\item \underline{\textit{Next event:}} degradation of machine $j$. The corresponding transition rate is
$$q^{\pmb{a}}\big((\pmb{F}, \pmb{P}, \pmb{C}, j), (\pmb{F}-\pmb{e}_x\mathbb{1}\{y=-1\}-\pmb{e}_y, \pmb{P}+\pmb{e}_x, \pmb{C}^{j, N-1}, j)\big) = (1-\alpha_N)\mu_N.$$
\item \underline{\textit{Next event:}} repair of machine $l$, $C_l=0$. The corresponding transition rate is
$$q^{\pmb{a}}\big((\pmb{F}, \pmb{P}, \pmb{C}, j\big), (\pmb{F}-\pmb{e}_x\mathbb{1}\{y=-1\}-\pmb{e}_y, \pmb{P}+\pmb{e}_x, {\pmb{C}^{j, N}}^{l, N}, l)\big) = \mu_0.$$
\end{enumerate}

\noindent
\textbf{Type 3} \underline{\textit{Last event:}} degradation or repair of machine $j \in \mathcal{J}$; \underline{\textit{action:}} $\pmb{a}=(x, y, z) \in \mathcal{A}_2(\pmb{F}, \pmb{P}, \pmb{C}, j)$.\\
The last event is either a repair of machine $j$ or a degradation that is not a failure. No preventive maintenance is done, so $x=-1$. However, a relocation of one spare part between any two warehouses $y$ and $z$ is possible. The state of the system immediately after the action is taken is $\pmb{X} = (\pmb{F}-\pmb{e}_y+\pmb{e}_z, \pmb{P}, \pmb{C}, j)$. The transition rates are defined as follows.
\begin{enumerate}
\item \underline{\textit{Next event:}} replenishment at $k$, $P_k>0$. The corresponding transition rate is
$$q^{\pmb{a}}\big((\pmb{F}, \pmb{P}, \pmb{C}, j), (\pmb{F}-\pmb{e}_y+\pmb{e}_z+\pmb{e}_k, \pmb{P}-\pmb{e}_k, \pmb{C}, 0)\big) = P_k\gamma.$$
\item \underline{\textit{Next event:}} failure of machine $l$, $C_l>0$. The corresponding transition rate is
$$q^{\pmb{a}}\big((\pmb{F}, \pmb{P}, \pmb{C}, j), (\pmb{F}-\pmb{e}_y+\pmb{e}_z, \pmb{P}, \pmb{C}^{l, 0}, l)\big) = \alpha_{C_l}\mu_{C_l}.$$
\item \underline{\textit{Next event:}} degradation of machine $l$, $C_l>1$. The corresponding transition rate is
$$q^{\pmb{a}}\big((\pmb{F}, \pmb{P}, \pmb{C}, j), (\pmb{F}-\pmb{e}_y+\pmb{e}_z, \pmb{P}, \pmb{C}^{l, C_l-1}, l)\big) = (1-\alpha_{C_l})\mu_{C_l}.$$
\item \underline{\textit{Next event:}} repair of machine $l$, $C_l=0$. The corresponding transition rate is
$$q^{\pmb{a}}\big((\pmb{F}, \pmb{P}, \pmb{C}, j), (\pmb{F}-\pmb{e}_y+\pmb{e}_z, \pmb{P}, \pmb{C}^{l, N}, l)\big) = \mu_0.$$
\end{enumerate}

\noindent
\textbf{Uniformization}

To be able to compute the optimal policy, we uniformize our Markov process $\boldsymbol X^{\boldsymbol a}(t)$. To do so, we introduce the constant $\tau = \gamma K + J \max_{n \in \{0, ... ,N\}}{\mu_n}$ that is larger than the total transition rate from any state, and add the following dummy transitions that make the total transition rate from any state equal to $\tau$.
We add the following dummy transitions for each type of transition described above with $\pmb{a}=\{-1, -1, -1\}$.\\
\noindent
Type 1:
\begin{equation*}
\begin{split}
q^{\pmb{a}}\big((\pmb{F}, \pmb{P}, \pmb{C}, 0), (\pmb{F}, \pmb{P}, \pmb{C}, 0)\big) &= \tau - \sum_{k=1}^{I}P_k\gamma-\sum_{l \in \mathcal{J}_1}\alpha_{C_l}\mu_{C_l} - \sum_{l \in \mathcal{J}_2}(1-\alpha_{C_l})\mu_{C_l} - \mu_0 |\mathcal{J}_0|\\
&=\tau - \sum_{k=1}^{I}P_k\gamma-\sum_{l \in \mathcal{J}_1}\mu_{C_l} - \mu_0 |\mathcal{J}_0|\\
&=\tau - \sum_{k=1}^{I}P_k\gamma-\sum_{l \in \mathcal{J}}\mu_{C_l},
\end{split}
\end{equation*}
where $\mathcal{J}_0 = \{l \in \mathcal{J}: C_l = 0\}$, $\mathcal{J}_1 = \{l \in \mathcal{J}: C_l > 0\}$ and $\mathcal{J}_2 = \{l \in \mathcal{J}: C_l > 1\}$.\\
\noindent
Type 2.1:
\begin{equation*}
\begin{split}
q^{\pmb{a}}\big((\pmb{F}, \pmb{P}, \pmb{C}, j), (\pmb{F}-\pmb{e}_x\mathbb{1}\{y=-1\}-\pmb{e}_y, \pmb{P}+\pmb{e}_x, \pmb{C}, 0)\big) = \tau - \sum_{k=1}^{I}(\pmb{P}+\pmb{e}_x)_k\gamma-\sum_{l \in \mathcal{J}}\mu_{C_l}.
\end{split}
\end{equation*}
\noindent
Type 2.2:
\begin{equation*}
\begin{split}
q^{\pmb{a}}\big((\pmb{F}, \pmb{P}, \pmb{C}, j), (\pmb{F}-\pmb{e}_x\mathbb{1}\{y=-1\}-\pmb{e}_y, \pmb{P}+\pmb{e}_x, \pmb{C}^{j, N}, 0)\big) = &\tau - \sum_{k=1}^{I}(\pmb{P}+\pmb{e}_x)_k\gamma\\
&-\sum_{l \in \mathcal{J}}\mu_{C_l}+\mu_{C_j}-\mu_N.
\end{split}
\end{equation*}
\noindent
Type 3:
\begin{equation*}
\begin{split}
q^{\pmb{a}}\big((\pmb{F}, \pmb{P}, \pmb{C}, j), (\pmb{F}-\pmb{e}_y+\pmb{e}_z, \pmb{P}, \pmb{C}, 0)\big) = \tau - \sum_{k=1}^{I}P_k\gamma-\sum_{l \in \mathcal{J}}\mu_{C_l}.
\end{split}
\end{equation*}

\subsection{Costs}\label{sec:model:costs}

We consider a generic cost structure that would allow to study various types of settings and examine the effects of different actions on the optimal policy. We incorporate the following cost components that are common in research literature on lateral transshipments and CMB. Let $c_{cs}$ and $c_{ps}$ be fixed setup costs for corrective and preventive maintenance, respectively, given that a spare part is dispatched from a local warehouse. In case a spare part is dispatched from the central warehouse, the costs $c_e$ are incurred, independent of the type of maintenance. Let $c_{rs}$ denote the setup costs incurred per relocation, and $c_r$ - the replenishment setup costs. Assume that for each pair of local warehouse $i$ and machine $j$ the corresponding response time $R_{ij}$ is deterministic and known. A fixed penalty $c_{cl}$ is incurred if response time to a failed machine is larger than a given time threshold $t^*$, and an extra penalty of $c_{cp}$ per time unit of delay over $t^*$. We assume that, in case a spare part is dispatched from the central warehouse, response time is always smaller than the time threshold $t^*$, independent of the machine. The immediate costs of action $\boldsymbol a(\boldsymbol X) = (x, y, z)$ in state $\boldsymbol X = (\boldsymbol F, \boldsymbol P, \boldsymbol C, j)$ can be computed as follows:
\begin{equation}\label{eq:mdp_costs}
\mbox{\footnotesize$\displaystyle
c(\boldsymbol X,\boldsymbol a(\boldsymbol X)) = \left\{
\begin{array}{ll}
c_e		&	{\rm if~} x=0,\\
c_{cs}+c_r+(c_{cl}+c_{cp}(R_{xj}-t^*))\mathbb{1}\{R_{xj}>t^*\}+c_{rs}\mathbb{1}\{y>0\} &	{\rm if~} x>0 \text{ and } C_j=0,\\
c_{ps}+c_r+c_{rs}\mathbb{1}\{y>0\} &	{\rm if~} x>0 \text{ and } C_j>0,\\
c_{rs}\mathbb{1}\{y>0\}			&	{\rm if~} x=-1.
\end{array}\right.
$}
\end{equation}

\subsection{Optimality Equations}

We formulate the problem as an infinite-horizon discounted MDP. Let $V(\pmb{X})$ denote the expected total discounted costs under the optimal policy, when starting in state $\pmb{X}$. Then $V(\pmb{X})$ satisfies the Bellman equations:
\begin{equation}\label{eq:mdp_bellman}
V(\pmb{X}) = \sup_{\pmb{a}\in \mathcal{A}(\pmb{X})}\left\{c(\pmb{X}, \pmb{a}) + \sum_{\pmb{X'}\in \mathcal{S}}\lambda p(\pmb{X'}|\pmb{X}, \pmb{a})V(\pmb{X'}) \right\},
\end{equation}
where $\lambda<1$ is a discount factor.

\section{Numerical Experiments}\label{sec:num}
In this section we conduct a number of experiments to study the performance and the structure of the optimal policy. To compute the optimal policy, we use the policy iteration algorithm with the maximum number of iterations set to 1000. All experiments are run in Python 3.7 on a computer with 8 GB RAM, Intel Core i5-5250U 1.6 GHz processor, running Linux Fedora 30.

\subsection{Experimental Setup}\label{sec:num:setup}
\textit{Parameters.} The following parameters are fixed throughout all experiments: the time threshold $t^*=10$, the discount factor $\lambda=0.95$, the number of warehouses $I=2$, the number of machines $J=2$, and the inventory level $K=2$. Note that we only consider small problem instances, as due to the curse of dimensionality, it would be infeasible to derive optimal policy for multiple instances and for a wide range of parameter settings. We also assume $\mu_i=1~(i = 0, 1,...,N)$, and $\alpha_i = 0~(i = 2, ..., N)$. To study the system performance under the different levels of workload, we introduce the load parameter $\rho = \frac{J}{N \gamma K}$. For given values of $\rho$ and $N$, we adjust $\gamma$ accordingly.

\begin{figure}
\begin{subfigure}{.49\textwidth}
\centering
\includegraphics[width=\textwidth]{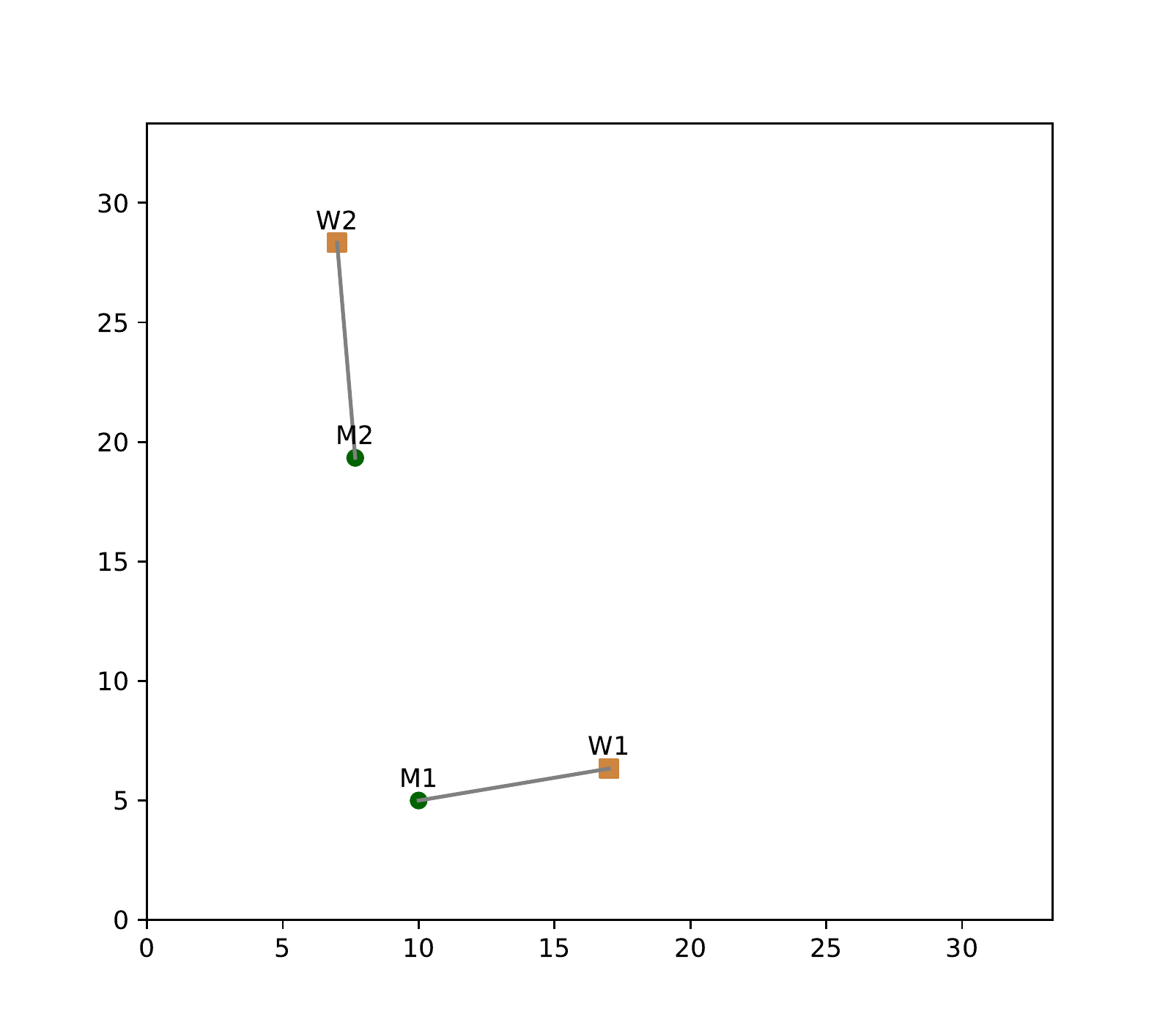}
\caption{}
\label{fig:map_1}
\end{subfigure}\hfill
\begin{subfigure}{.49\textwidth}
\centering
\includegraphics[width=\textwidth]{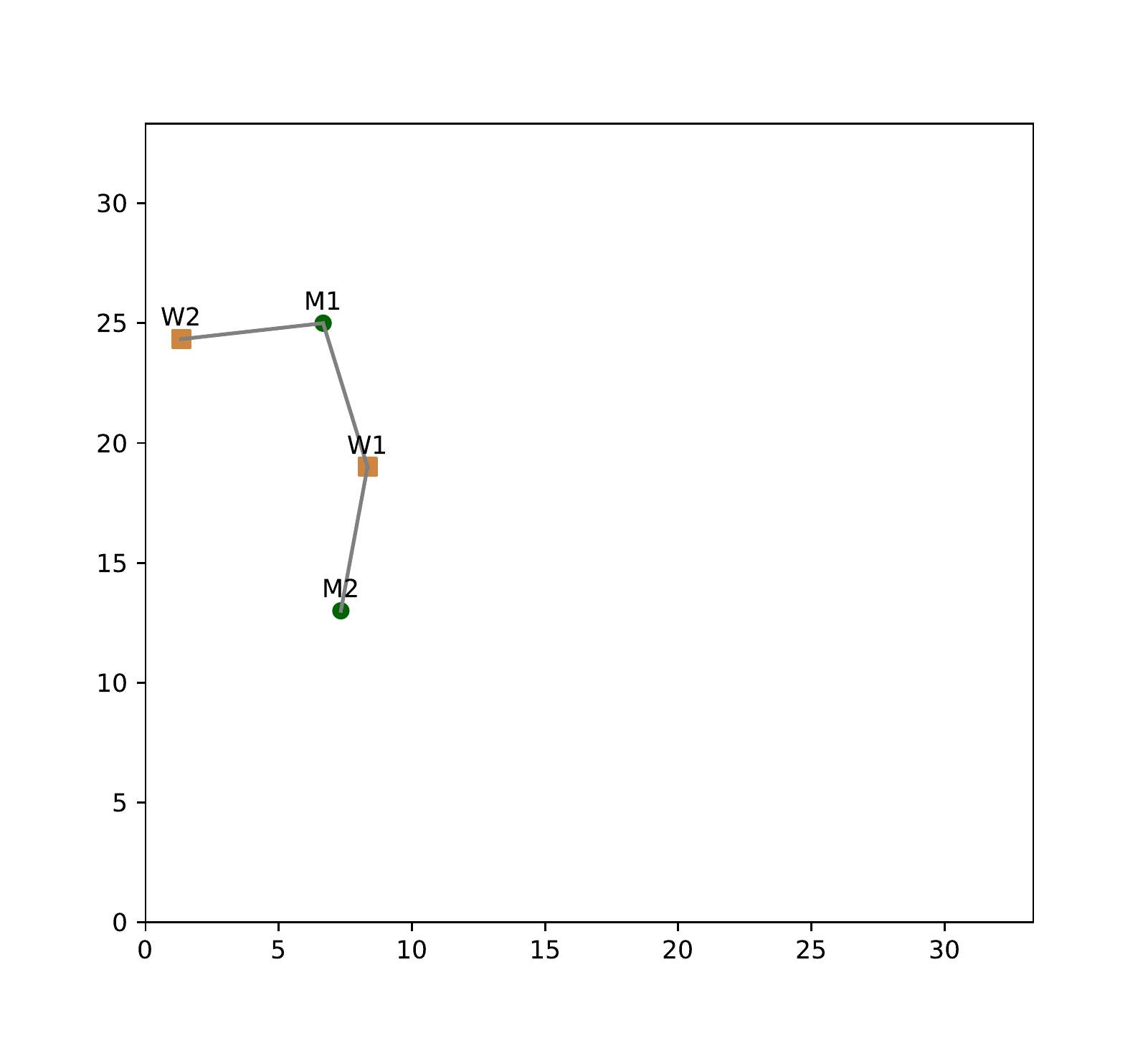}
\caption{}
\label{fig:map_2}
\end{subfigure}
	\caption{Examples of problem instances. The warehouses are connected to the machines that are reachable within $t^*$ time units}
    \label{fig:maps}
\end{figure}

\textit{Response times.} An important component of a problem instance is a matrix $\pmb{R}$ of fixed response times $R_{ij}$ between each pair of warehouse $i$ and machine $j$. The matrix $\pmb{R}$ is used in equation~\eqref{eq:mdp_costs} to compute the immediate costs. For a given random seed, we construct it as follows. Machines and warehouses are allocated at random within a square of size $33\times 33$ in terms of time units, such that each warehouse is within $t^*=10$ time units from at least one machine, and each machine is within $t^*=10$ time units from at least one warehouse. The response times $R_{ij}$ are then computed as the corresponding Euclidean distances. Figure~\ref{fig:maps} presents two examples of problem instances used in this study.

\textit{Policy types.} To study the effects of different types of actions on the policy performance, and in particular, the effects of condition-based maintenance, we introduce the following five types of policies that are defined by limiting the original action spaces~\eqref{eq:mdp_actions1} and~\eqref{eq:mdp_actions2} as follows:
\begin{enumerate}
\item \textbf{Closest-First corrective maintenance (CF).} With this policy type only corrective maintenance is done, using the closest available spare part in terms of response time. The central warehouse stock is used only if all local warehouses are empty. The corresponding action space is defined as follows:
\begin{flalign*}
\mbox{\small$\displaystyle
\mathcal{A}_{CF}(\boldsymbol X) =
\begin{cases}
\{(x, y, z)|x = \argmin_{i\in \mathcal{W}(\pmb{X})}R_{ij}, \ y =-1, \ z = -1\}, & C_j = 0, \ \mathcal{W}(\pmb{X}) \neq \emptyset, \\
\{(0, -1, -1)\}, & C_j = 0, \ \mathcal{W}(\pmb{X}) = \emptyset, \\
\{(-1, -1, -1)\}, & \text{otherwise}.
\end{cases}
$}
\end{flalign*}%
Note that with this policy type there is exactly one action per state. So, there is no need to use policy iteration, and the value function can be obtained by solving a set of linear equations:
\begin{equation*}
V(\pmb{X}) = c(\pmb{X}, \pmb{a}_{CF}(\pmb{X})) + \sum_{\pmb{X'}\in \mathcal{S}}\lambda p(\pmb{X'}|\pmb{X}, \pmb{a}_{CF}(\pmb{X}))V(\pmb{X'}),
\end{equation*}
where $\pmb{a}_{CF}(\pmb{X})$ is the action taken in state $\pmb{X}$ under the CF policy.

\item \textbf{Optimal Corrective maintenance (OC).} With this policy type only corrective maintenance is done that is subject to optimization as in equations~\eqref{eq:mdp_bellman}. The corresponding action space is defined as follows:
\begin{flalign*}
\mbox{\small$\displaystyle
\mathcal{A}_{OC}(\boldsymbol X) =
\begin{cases}
\{(x, y, z) \in \mathcal{A}_1(\boldsymbol X)| y=-1\}, & C_j = 0, \\
\{(-1, -1, -1)\}, & \text{otherwise}.
\end{cases}
$}
\end{flalign*}%

\item \textbf{Optimal Corrective maintenance with Relocation (OCR).} With this policy both corrective maintenance and relocation actions are optimized, given that no preventive maintenance is done. Relocation is also allowed upon a change of a machine state that is not a failure. The corresponding action space is defined as follows:
\begin{flalign*}
\mbox{\small$\displaystyle
\mathcal{A}_{OCR}(\boldsymbol X) =
\begin{cases}
\mathcal{A}_1(\boldsymbol X), & C_j = 0, \\
\mathcal{A}_2(\boldsymbol X), & \text{otherwise}.
\end{cases}
$}
\end{flalign*}%

\item \textbf{Optimal Corrective and Preventive maintenance (OCP).} With this policy type both corrective and preventive maintenance are subject to optimization, and no relocations are allowed. The corresponding action space is defined as follows:
\begin{flalign*}
\mbox{\small$\displaystyle
\mathcal{A}_{OC}(\boldsymbol X) =
\{(x, y, z) \in \mathcal{A}_1(\boldsymbol X)| y=-1\}
$}
\end{flalign*}%

\item \textbf{Optimal Corrective and Preventive maintenance with Relocation (OCPR).} The last policy type corresponds to the full action space defined in Section~\ref{sec:model:actions}:
\begin{flalign*}
\mbox{\small$\displaystyle
\mathcal{A}(\boldsymbol X) = \mathcal{A}_1(\boldsymbol X, j) \cup \mathcal{A}_2(\boldsymbol X, j).
$}
\end{flalign*}%
\end{enumerate}

\textit{Performance measures.} Solving the Bellman equations gives a value function $\pmb{V}$ with the total expected discounted costs per state. To measure the policy performance we use the weighted average of the components of $\pmb{V}$, where the steady state probabilities under the optimal policy are used as the weights. The steady state probabilities vector $\pmb{\pi}$ is computed by solving the system of linear equations: 
\begin{flalign*}
\begin{cases}
\pmb{\pi}\pmb{P} = \pmb{\pi},\\
\sum_{i=1}^{|\mathcal{S}|}\pmb{\pi}=1,
\end{cases}
\end{flalign*}%
where $\pmb{P}$ is the matrix of transition probabilities under the optimal policy. The weighted average of the value function is denoted by $\upsilon = \pmb{\pi}'\pmb{V}$.

For the policy types 2 to 5 we also report the relative improvement over the CF policy denoted by $\Delta$. For example, for the optimal OCPR policy we define
$$\Delta_{OCPR} = \frac{\upsilon_{CF} - \upsilon_{OCPR}}{\upsilon_{CF}}\times 100\%.$$

\subsection{Different Cost Settings}

In this section we study the performance of the different policies depending on the cost setting and the load. The general assumption we make when choosing the values for cost parameters are aligned with the research literature on spare parts management and are as follows. The setup costs for relocation and dispatching preventively are lower than the setup costs for dispatching correctively from a local warehouse. Dispatching from the central warehouse has higher setup costs than corrective dispatching from a local warehouse. as it is supposed to be done only in emergency situations.
We study the policies' performance under different values of the load parameter $\rho$. For each combination of the cost parameters and the load $\rho$, 30 random instances are generated as described in Section~\ref{sec:num:setup}, and the average performance is computed. Table~\ref{tbl:cost_setups} presents the obtained results.

The immediate costs~\eqref{eq:mdp_costs} depend on seven parameters, and it is infeasible to cover all possible cases. Hence, we choose the following three different cost settings:
\begin{enumerate}
\item $C_p = 0.05$, $C_r = 0$, $C_{cs} = 1$, $C_{cl} = 1$, $C_{ps} = 0.2$, $C_{rs} = 0.2$, $C_e = 10$;
\item $C_p = 0.1$, $C_r = 0$, $C_{cs} = 10$, $C_{cl} = 1$, $C_{ps} = 0.2$, $C_{rs} = 0.2$, $C_e = 100$;
\item $C_p = 0$, $C_r = 0$, $C_{cs} = 0$, $C_{cl} = 1$, $C_{ps} = 0$, $C_{rs} = 0$, $C_e = 10$.
\end{enumerate}

Setting 1 corresponds to machinery of \textit{moderate criticality}, where it is important to address breakdowns within the given time limit. The delay in response time is also penalized, although not significantly. Setting 2 corresponds to \textit{critical machinery}, where breakdowns are very costly independent of response time. There is also a larger penalty for the delay in response time. Setting 3 corresponds to the case where breakdowns are \textit{not critical} as long as they are taken care of within the time limit. For all three cost settings we choose relocation and preventive maintenance setup costs to be noticeably lower than the corrective setup costs and equal to each other.

\begingroup
\setlength{\tabcolsep}{6pt} 
\renewcommand{\arraystretch}{1} 
\begin{table}[]
\small
\caption{Average performance of the policies per cost setting and load over 30 problem instances}
\label{tbl:cost_setups}
\centering
\begin{tabular}{|c|c|c|c|c|c|c|c|c|c|c|}
\hline
\multicolumn{2}{|c|}{}         & \multicolumn{5}{c|}{\textbf{Average} $\pmb{\upsilon}$}                                         & \multicolumn{4}{c|}{\textbf{Average} $\pmb{\Delta}$}                       \\ \hline
\textbf{Cost setting}        & $\pmb{\rho}$ & \textbf{CF} & \textbf{OC} & \textbf{OCR} & \textbf{OCP} & \textbf{OCPR} & \textbf{OC} & \textbf{OCR} & \textbf{OCP} & \textbf{OCPR} \\ \hline
\multirow{4}{*}{\textbf{1}} & \textbf{1}   & 7.19        & 7.12        & 6.57         & 6.98         & 6.57          & 1.0\%       & 8.6\%        & 2.9\%        & 8.7\%         \\ \cline{2-11} 
                            & \textbf{0.7} & 5.37        & 5.29        & 4.79         & 4.89         & 4.54          & 1.5\%       & 10.8\%       & 9.0\%        & 15.6\%        \\ \cline{2-11} 
                            & \textbf{0.5} & 4.02        & 3.94        & 3.50         & 3.48         & 3.16          & 2.1\%       & 13.0\%       & 13.6\%       & 21.5\%        \\ \cline{2-11} 
                            & \textbf{0.3} & 2.54        & 2.46        & 2.13         & 2.09         & 1.85          & 3.3\%       & 16.3\%       & 17.6\%       & 27.3\%        \\ \hline
\multirow{4}{*}{\textbf{2}} & \textbf{1}   & 63.67       & 63.62       & 63.19        & 62.18        & 61.89         & 0.1\%       & 0.8\%        & 2.3\%        & 2.8\%         \\ \cline{2-11} 
                            & \textbf{0.7} & 46.03       & 45.97       & 45.58        & 41.49        & 41.22         & 0.1\%       & 1.0\%        & 9.9\%        & 10.5\%        \\ \cline{2-11} 
                            & \textbf{0.5} & 33.30       & 33.24       & 32.90        & 28.10        & 27.86         & 0.2\%       & 1.2\%        & 15.6\%       & 16.3\%        \\ \cline{2-11} 
                            & \textbf{0.3} & 19.91       & 19.85       & 19.60        & 15.74        & 15.57         & 0.3\%       & 1.6\%        & 21.0\%       & 21.8\%        \\ \hline
\multirow{4}{*}{\textbf{3}} & \textbf{1}   & 5.24        & 5.15        & 4.48         & 5.15         & 4.48          & 1.7\%       & 14.4\%       & 1.7\%        & 14.4\%        \\ \cline{2-11} 
                            & \textbf{0.7} & 3.55        & 3.45        & 2.84         & 3.28         & 2.84          & 2.8\%       & 20.0\%       & 7.6\%        & 19.8\%        \\ \cline{2-11} 
                            & \textbf{0.5} & 2.38        & 2.28        & 1.74         & 2.01         & 1.62          & 4.5\%       & 26.8\%       & 15.7\%       & 32.0\%        \\ \cline{2-11} 
                            & \textbf{0.3} & 1.25        & 1.15        & 0.74         & 0.93         & 0.63          & 8.3\%       & 40.6\%       & 25.7\%       & 49.8\%        \\ \hline
\end{tabular}
\end{table}
\endgroup

For each cost setup and each value of the load parameter $\rho$ we generate 30 different problem instances. For  instances and compute the average performance of the optimal policy for each of the policy types. The obtained results are reported in Table~\ref{tbl:cost_setups}. We observe that the largest improvement over the CF policy is obtained under the cost setting 3 for all of the other policy types. Note that for the cost setting 2 optimal corrective maintenance and relocation have only marginal effect, while preventive maintenance results in a significant reduction in costs. For cost settings 1 and 3 doing relocations (OCR) has a bigger effect than doing preventive maintenance (OCP).

\subsection{Importance of Better Condition Diagnostics}
With better diagnostics we can more accurately identify at which point of a degradation process a machine is. We model improvement in diagnostics by decomposing the degradation process in a larger number of intermediate steps, that is, by increasing $N$ while keeping the load $\rho$ fixed. We consider the cost setting 1, and as before, use 30 instances per parameter setting. Table~\ref{tbl:varying_N} shows that the average $\Delta$ of OCP and OCRP policies increases significantly with $N$ for different loads. This means the contribution of preventive maintenance grows with $N$, demonstrating the importance of accurate diagnostics.

\begingroup
\setlength{\tabcolsep}{6pt} 
\renewcommand{\arraystretch}{1} 
\begin{table}[]
\small
\caption{Average performance of the policies for different $N$ over 30 problem instances}
\label{tbl:varying_N}
\centering
\begin{tabular}{|c|c|c|c|c|c|c|c|c|c|c|}
\hline
\multicolumn{2}{|c|}{}                     & \multicolumn{5}{c|}{\textbf{Average} $\pmb{\upsilon}$}                                         & \multicolumn{4}{c|}{\textbf{Average} $\pmb{\Delta}$}                       \\ \hline
$\pmb{\rho}$                  & \textbf{N} & \textbf{CF} & \textbf{OC} & \textbf{OCR} & \textbf{OCP} & \textbf{OCPR} & \textbf{OC} & \textbf{OCR} & \textbf{OCP} & \textbf{OCPR} \\ \hline
\multirow{5}{*}{\textbf{1}}   & \textbf{2} & 7.19        & 7.12        & 6.57         & 6.98         & 6.57          & 1.0\%       & 8.6\%        & 2.9\%        & 8.7\%         \\ \cline{2-11} 
                              & \textbf{3} & 6.69        & 6.64        & 6.20         & 5.71         & 5.44          & 0.8\%       & 7.4\%        & 14.7\%       & 18.7\%        \\ \cline{2-11} 
                              & \textbf{4} & 6.03        & 5.99        & 5.62         & 4.70         & 4.55          & 0.7\%       & 6.8\%        & 22.1\%       & 24.5\%        \\ \cline{2-11} 
                              & \textbf{5} & 5.42        & 5.39        & 5.07         & 3.73         & 3.65          & 0.6\%       & 6.4\%        & 31.2\%       & 32.7\%        \\ \cline{2-11} 
                              & \textbf{6} & 4.89        & 4.86        & 4.59         & 3.10         & 3.23          & 0.6\%       & 6.2\%        & 36.6\%       & 33.9\%        \\ \hline
\multirow{5}{*}{\textbf{0.5}} & \textbf{2} & 4.02        & 3.94        & 3.50         & 3.48         & 3.16          & 2.1\%       & 13.0\%       & 13.6\%       & 21.5\%        \\ \cline{2-11} 
                              & \textbf{3} & 4.06        & 4.00        & 3.95         & 2.98         & 2.93          & 1.6\%       & 2.8\%        & 26.6\%       & 28.0\%        \\ \cline{2-11} 
                              & \textbf{4} & 3.88        & 3.83        & 3.84         & 2.27         & 2.14          & 1.4\%       & 0.9\%        & 41.6\%       & 45.0\%        \\ \cline{2-11} 
                              & \textbf{5} & 3.64        & 3.59        & 3.39         & 1.74         & 1.75          & 1.2\%       & 6.8\%        & 52.0\%       & 51.7\%        \\ \cline{2-11} 
                              & \textbf{6} & 3.39        & 3.35        & 3.09         & 1.35         & 1.36          & 1.1\%       & 8.9\%        & 60.1\%       & 59.9\%        \\ \hline
\end{tabular}
\end{table}
\endgroup

\subsection{Balancing Relocation and Preventive Maintenance}
\begin{figure}
\begin{subfigure}{.5\textwidth}
\centering
\includegraphics[width=\textwidth]{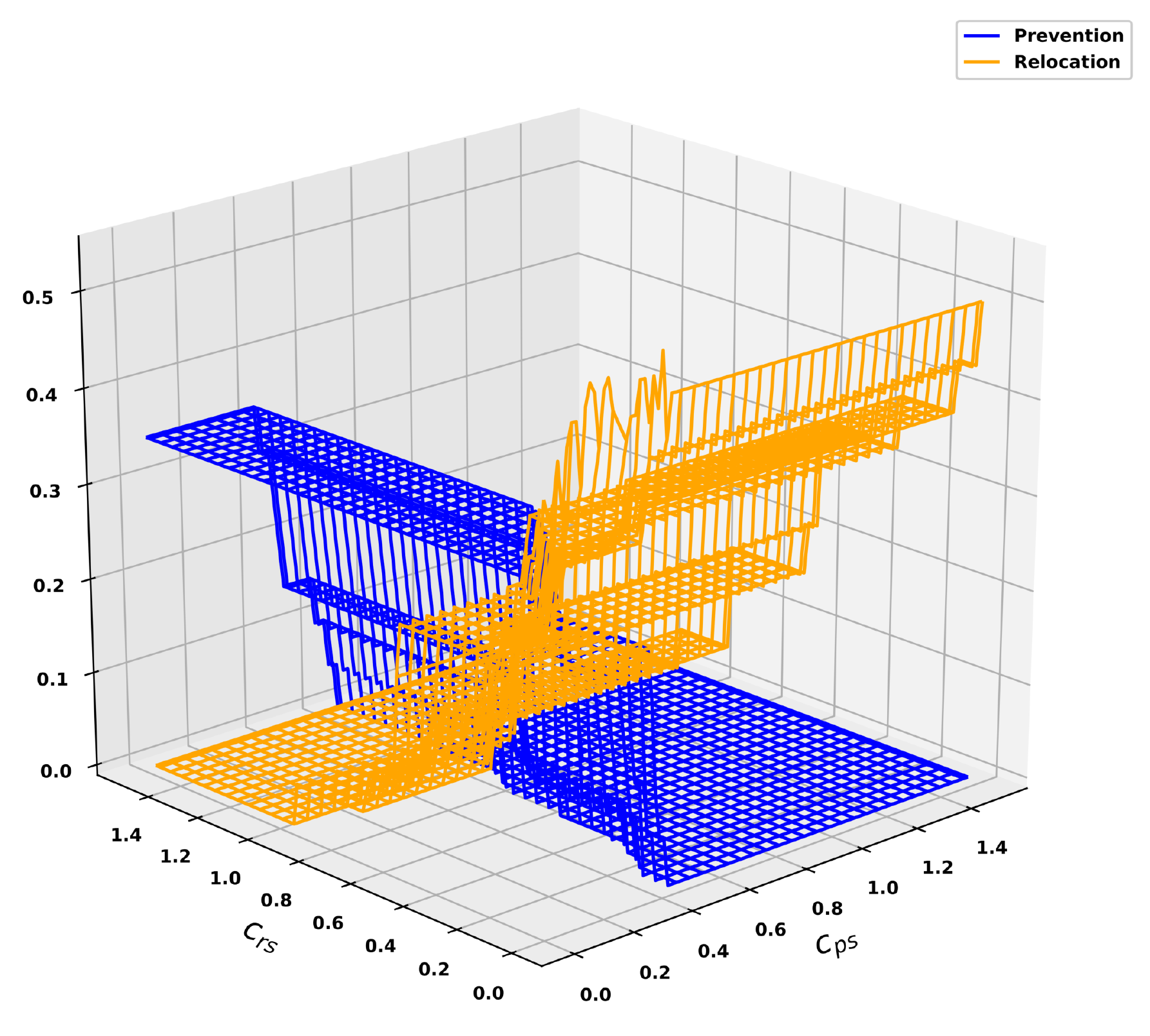}
\end{subfigure}\hfill
\begin{subfigure}{.5\textwidth}
\centering
\includegraphics[width=\textwidth]{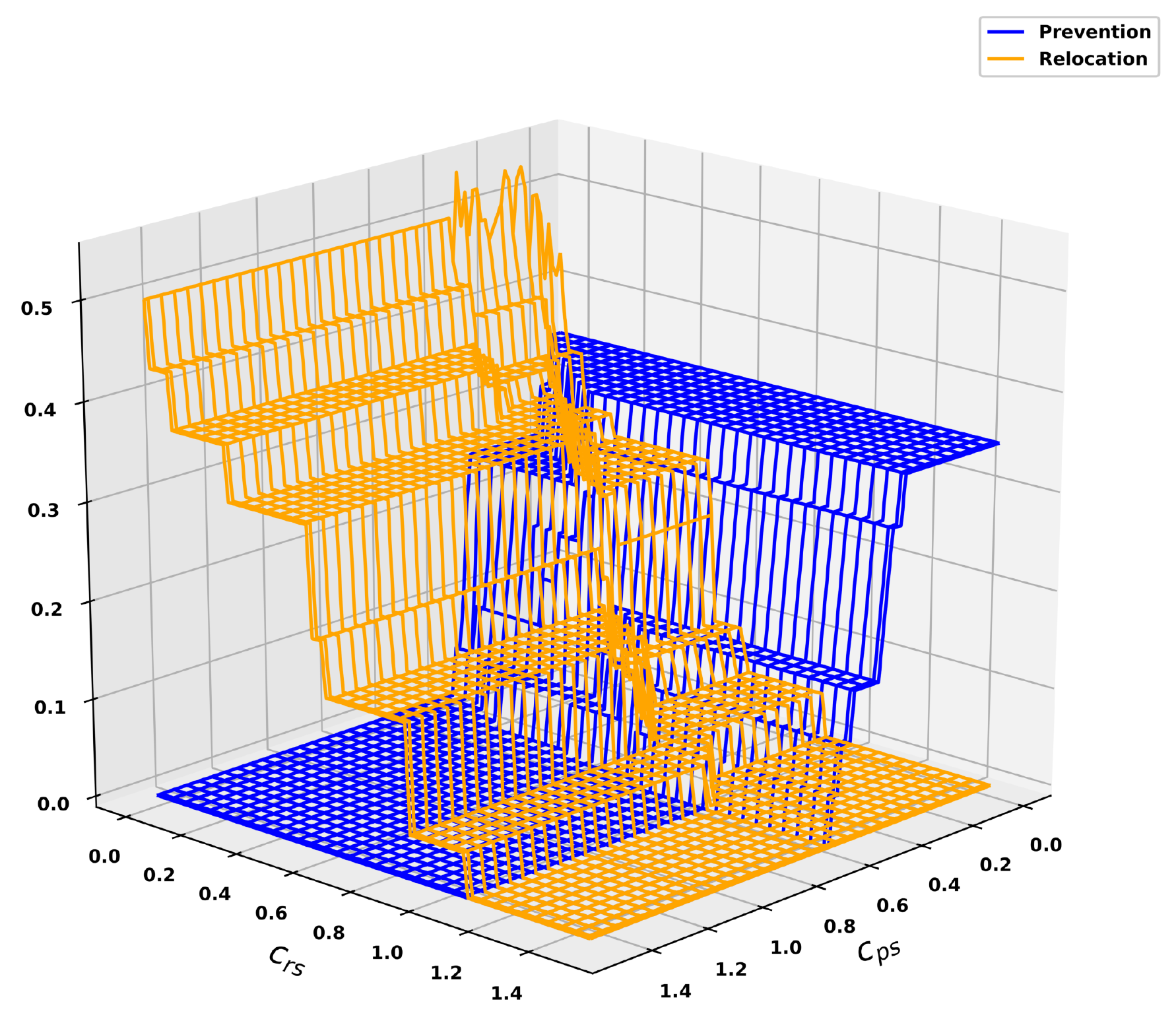}
\end{subfigure}
\caption{Relative number of states with relocation / preventive maintenance actions in the optimal OCPR policy as a function of $c_{ps}$ and $c_{rs}$}
\label{fig:prev_vs_reloc}
\end{figure}

In this section we show an example of how relocation actions are balanced with preventive maintenance actions in the optimal OCPR policy. We consider the problem instance from Figure~\ref{fig:map_1} and fix the parameters $\rho=0.5$ and $N=2$. We vary the cost components $c_{ps}$ and $c_{rs}$ in range $[0, 1.5]$ each, with other components fixed as in the cost setting 1. For each combination we compute the total number of states where relocation (preventive maintenance) is done in the optimal policy, divided by the total number of states where relocation (preventive maintenance) is possible. In Figure~\ref{fig:prev_vs_reloc} this metric is plotted against $c_{ps}$ and $c_{rs}$ for both prevention and relocation actions. We observe that both types of actions take place in the optimal policy while both $c_{ps}$ and $c_{rs}$ are relatively low. When one of the two cost components increases, the optimal policy leans towards either of the two with lower setup costs, and when both 
$c_{ps}$ and $c_{rs}$ are large, the optimal policy does not include either relocation or preventive maintenance actions.

\section{Conclusion}\label{sec:conclusion}

The work in this paper is a pioneering contribution to the field of dynamic spare parts management. We  introduce the concept of condition-based maintenance into the problem of dynamic dispatching and relocation of spare parts in a service network, and study the effects of this on the optimal policy. With the degradation process explicitly incorporated into the model, preventive maintenance of the machines and proactive relocation of spare parts become possible based on the current condition of all machines in the network as well as the availability and spatial distribution of resources.

We formulate the problem as an MDP, and study the optimal performance of various types of policies to evaluate the relative contribution of introducing CBM in a spare parts network. To that end, we conduct numerical experiments with a different cost settings, and show that the policies that use the information about the condition of the machines outperform those not doing so. We also demonstrate that better condition diagnostics can further improve the CBM based policy performance.

Due to the curse of dimentionality, solving MDP is computationally infeasible for large networks. Hence, in this work we only consider small problem instances. Given the benefits of introducing CBM on a network that we show in this paper, further research should focus on developing scalable heuristic approaches to the problem that would work for the problem instances of realistic sizes. Another interesting direction for further research is the parametric study of the degradation process. One could consider the effects of the corresponding parameters on the policy structure and its performance.

\noindent
\textbf{Acknowledgements}

This research was funded by an NWO grant, under contract number 438-15-506.

\bibliographystyle{plain}
\bibliography{references}

\appendix


\end{document}